\documentclass{article}
\title{On the two differences \\ \ $l_R(I^*/R)-l_R(R/I)$ \ and \
$rl_R(R/I)-l_R(I^*/R)$}

      \author{Anna Oneto$^{a,}${\footnote{Corresponding author.
{\em E-mail addresses}:  oneto@dimet.unige.it 
(A.Oneto),\qquad\qquad\qquad 
zatini@dima.unige.it (E.Zatini)}}
   \ and \  Elsa Zatini$^b$\\
{\small $^{a }$ {\it        Ditpem,
   Universit\`a di Genova,   P.le Kennedy, Pad. D 
-  I 16129 Genova (Italy);}} \\ {\small {\it$^b$
     Dima,  Universit\`a di Genova,  Via 
Dodecaneso 35 -  16146 Genova (Italy) }}}
\date{}

\expandafter\chardef\csname pre amssym.def  at\endcsname=\the\catcode`\@

\catcode`\@=11

\def\undefine#1{\let#1\undefined}
\def\newsymbol#1#2#3#4#5{\let\next@\relax
    \ifnum#2=\@ne\let\next@\msafam@\else
    \ifnum#2=\tw@\let\next@\msbfam@\fi\fi
    \mathchardef#1="#3\next@#4#5}
\def\mathhexbox@#1#2#3{\relax
    \ifmmode\mathpalette{}{\m@th\mathchar"#1#2#3}%
    \else\leavevmode\hbox{$\m@th\mathchar"#1#2#3$}\fi}
\def\hexnumber@#1{\ifcase#1 0\or 1\or 2\or 3\or 4\or 5\or 6\or 7\or 8\or
    9\or A\or B\or C\or D\or E\or F\fi}

\font\tenmsa=msam10
\font\sevenmsa=msam7
\font\fivemsa=msam5
\newfam\msafam
\textfont\msafam=\tenmsa
\scriptfont\msafam=\sevenmsa
\scriptscriptfont\msafam=\fivemsa
\edef\msafam@{\hexnumber@\msafam}

\font\tenmsb=msbm10
\font\sevenmsb=msbm7
\font\fivemsb=msbm5
\newfam\msbfam
\textfont\msbfam=\tenmsb
\scriptfont\msbfam=\sevenmsb
\scriptscriptfont\msbfam=\fivemsb
\edef\msbfam@{\hexnumber@\msbfam}

      \font\tengothic=eufm10
      \font\sevengothic=eufm7
      \newfam\gothicfam
      \textfont\gothicfam=\tengothic
      \scriptfont\gothicfam=\sevengothic
      \def\goth#1{{\fam\gothicfam #1}}
      \font\tenmsb=msbm10
      \font\sevenmsb=msbm7
      \newfam\msbfam
      \textfont\msbfam=\tenmsb
      \scriptfont\msbfam=\sevenmsb
      
\newtheorem{prop}{Proposition}[section]
\newtheorem{rem}[prop]{Remark}
\newtheorem{thm}[prop]{Theorem}
\newtheorem{coro}[prop]{Corollary}
\newtheorem{defin}[prop]{Definition}
\newtheorem{lemma}[prop]{Lemma}
\newtheorem{notat}[prop]{Notation}

\newtheorem{(*)}[prop]{}

\newcommand{\integ}{\mbox{${\sf Z\hspace{-1.3 mm}Z}$}}

\newcommand{\nat}{\mbox{${\rm I\hspace{-.6 mm}N}$}}
\newcommand{\Ib}{\mbox{$\overline{I}$}}

\newcommand{\R}{\mbox{$\overline{R}$}}

\newcommand{\ded}{\mbox{$\theta_D$}}

\newcommand{\m}{\mbox{${\goth m}$}}

\newcommand{\w}{\mbox{$\omega$}}

\newcommand{\I}{\mbox{$\  \Longrightarrow \ $}}
\newcommand{\II}{\mbox{$\  \Longleftrightarrow \ $}}

\newcommand{\LL}{\mbox{$\Lambda$}}

\font\tengothic=eufm10
\font\sevengothic=eufm7
\newfam\gothicfam
\textfont\gothicfam=\tengothic
\scriptfont\gothicfam=\sevengothic
\def\goth#1{{\fam\gothicfam #1}}

\begin{document}

\maketitle

 {\bf Abstract}.  Let  $R $    be a one-dimensional local Noetherian
      domain, which is supposed   analytically irreducible and residually
rational, and let $I$ be a proper ideal of $R$. Our purpose   is to study
the two numbers \vspace{-0.2cm}
  $$   a(I):=l_R(I^*/R)-l_R(R/I)\ \vspace{-0.2cm} $$
    $$b(I):=rl_R(R/I)-l_R(I^*/R) \vspace{-0.2cm}$$
already considered in the literature under various points of view. The
basic idea is the expression of these invariants
in terms of the type sequence.

  \section{Introduction.}

   Let ($R, \m$)   be a one-dimensional local Noetherian
      domain with
residue field $k$ and quotient field $K$, which
is analytically irreducible and residually
rational. We denote by:
\par
\vspace{ 0.1cm}
\  \ $\begin{array}{rcl}  \R   &{\rm the} &   normalization
\ {\rm  of} \ \ R,  \ \ \R  =k[[t]];\\ \omega &{\rm a} &     canonical\
module \
\
{\rm of }\ \  R\ \ {\rm such\ that \ \ }
\ R\subseteq
\w\subseteq
\R;\\    \gamma&:=& R:\R  \ \ \ {\rm the} \  conductor\ ideal \ \  {\rm of}
\ R\ {\rm in}   \ \R; \\
      c&:=& l_R(\R/\gamma ),  {\rm \ \ so \ that\ }
      \ \gamma   =t^c\overline R;\\
      \delta&:=& l_R(\R/R) \    {\rm  \ the\ }\ singularity\
degree \ {\rm of} \ R;\\
n&:=&c-\delta\ =\ l_R(R/\gamma );\\
r & :=& l_R(R:\m/R) \ \ {\rm the} \ \ Cohen-Macaulay\ type \ \  {\rm of} \ R;\\
I^* & := &R:_KI  \ \ {\rm the} \  dual \ {\rm of \ the \ fractional \
ideal} \ I;\\
\ded &:=& \w^* \ \   {\rm the} \  Dedekind \ different  \ {\rm of} \
R.\end{array}   $ \\ \vspace{-0.1cm}

\noindent Given any proper nonzero  ideal $I$ of $R$,
we use the notion
of {\it  type sequence}   (see \ref{ts}) \ in order to get
       informations about  the two numerical
differences:\vspace{-0.2cm}
$$  a(I):=l_R(I^*/R)-l_R(R/I)\ \ \vspace{-0.2cm}$$
    \centerline{$b(I):=rl_R(R/I)-l_R(I^*/R).$}
Having in mind the
     Gorenstein case  with the following well-known equivalent
characterizations \  (see \cite{b}; \cite{j};  \cite{ma}, Theorem 13.1):

    $R\  is\ Gorenstein\ \iff\ \w=R\ \iff \
r=1\ \iff   2\delta-c=0$  \par $ \ \ \ \ \ \ \ \ \ \ \ \ \ \ \ \ \ \ \ \ \
\ \ \  \iff   a(I)=0\ \
for\
every \ nonzero\ proper \ ideal\ I $     \\ \vspace{-0.2cm}

\noindent  we get  similar
characterizations     for    $almost\
Gorenstein$
rings
(see Theorem
\ref{Matlis13.1}):  \vspace{-0.2cm}
$$R \ is \ almost  \ Gorenstein \quad \iff\quad \m\w=\m\quad\iff
\ r-1=2 \delta -c \vspace{-0.2cm}$$   $ 
\qquad\iff a(I)=r-1 -  l_R(I^{**}/I) \ \ for\
every \ nonzero  \ proper
\ ideal\ I.$  \\ \vspace{-0.2cm}

\noindent In the general case  a direct
calculation gives immediately that 
\vspace{-0.2cm}$$a(I)\leq 2 \delta -c 
\vspace{-0.2cm}$$
for every   nonzero ideal
$I$.  The close relation with the type sequence
of $R$   \\
  \centerline{$ 2 \delta -c =  \sum_{h=1}^n \;(r_h-1)$}
  induces us to search which elements of this sequence contribute in our
invariants. This is discussed in   Section (4).
First, in Theorem \ref{r.nelleformule}, we obtain
the
formulas:   \vspace{-0.2cm}
$$(A)\hspace{2cm} a(I)   =\  \sum_{ h \in \{1,...,n\} \setminus V^I }\
(r_h-1) - \ l_R(I^{**}/I) \ - d(I) \hspace{0.7cm}  \vspace{-0.2cm}$$
$$(B)\hspace{2cm} b(I)= \sum_{ h \notin V^I}\ (r-r_h) + \  r
l_R(I^{**}/I) + d(I) \hspace{2cm}   \vspace{-0.2cm}  $$

\noindent where $V^I\subseteq  \nat$  is a subset in biunivocal  correspondence
with the valuations of the ideal $I^{**}$ and
$d(I)$  is a non-negative invariant
   (see \ref{def.d}),   closely related to the type sequence
$[r_1,...,r_n]$ of $R$ and to the valuations of
$I^{**}\! $. \\ Successively  lower and upper bounds and vanishing
conditions for the invariants $a(I)$ and $b(I)$
are derived directly from these expressions.
    For
instance  \vspace{-0.2cm}
$$(A_1) \hspace{2cm} a(I) \leq
(r-1)l_R\big ( R/ (I^{**}+\ded)\big) - \ 
l_R(I^{**}/I) \hspace{1cm} \vspace{-0.2cm}$$

\noindent which  improves   the inequality  \ $a(I) \leq (r-1) l_R(R/I)   $ \
obtained by
J\"ager in \cite{j}, Korollar 3, (2) \ \ and   \vspace{-0.2cm}
$$(A_2)\hspace{2cm} a(I)\geq
(r-1)  - \ l_R(I^{**}/I) \ - d(I) \hspace{2.5cm} \vspace{-0.2cm}$$

\noindent which gives a  sufficient
condition   for the positivity of $a(I)$. \\  We recall that
in \cite{b}, Anm.5,  R. Berger   conjectured that always \ $a(I)
\geq0$, \ but
there are counterexamples, we cite   the
following, exhibited  by  J\"ager in \cite{j}: \\ \vspace{-0.2cm}

\noindent if \
$R=k[[t^9,t^{15},t^{17},t^{23},t^{25},t^{29},t^{31}]]$
\  and \ $I= (t^{38},t^{44},t^{50})$, \ then
\
$a(I)=-1$. \\ \vspace{-0.2cm}

   \noindent From the preceding $(A_2)$ it turns out that  \  $a(I)
\geq r-1 \geq 0$   for every integrally
closed ideal  $I$,  because   this
condition implies that \  $I=I^{**}$ and also that
$d(I)=0$. \ The same holds \ for every
ideal $I$ such that
\
$\w \subseteq I:I$. \\ If $R$ is $ almost\ Gorenstein$, then $a(I)=2 \delta
-c\geq 0$ for every reflexive ideal $I$.\\ \vspace{-0.2cm}

Formula (B), by giving
$b(I)$ as a sum of non negative terms,   provides
the fact that always  \vspace{-0.1cm}$$b(I)\geq r \  l_R(I^{**}/I) \geq
0\vspace{-0.1cm}$$ and also the
following
vanishing condition:  \vspace{-0.2cm}
$$(VC) \hspace{1cm}  b(I)=0 \ \iff\  I^{**}=I, \ d(I)=0,\  \
r_i=r \ \ {\rm for \ all} \
i\notin V^I   \vspace{-0.2cm}$$

\noindent Since by definition \vspace{-0.1cm} 
$$a(I)+b(I)=(r-1)\ \! l_R(R/I) \vspace{-0.1cm} 
$$  it is
clear that inequalities $(A_1)$ and  $(A_2)$ may be read  respectively as a
lower and an upper bound   for $b(I)$. We explicit these for the
convenience of the reader. \vspace{-0.2cm}
$$(B_1) \hspace{2cm} b(I) \geq
(r-1)l_R( I^{**}\!\!+\ded /\ \! I )+ \ 
l_R(I^{**}/I)\hspace{1.8cm} \vspace{-0.2cm} $$ 
\vspace{-0.2cm}
$$(B_2) \hspace{2cm}  b(I) \leq (r-1)(l_R(R/I)-1) + l_R(I^{**}/I) +
d(I)\hspace{1cm}   \vspace{-0.1cm} $$

\noindent In the literature more attention has
been reserved to the particular
case $I=\gamma$. \ Notice that \vspace{-0.1cm} 
$$ a(\gamma)= 2\delta-c\quad { \rm and}
\quad b(\gamma)= r(c-\delta)-\delta \vspace{-0.1cm} $$ As concerns the number
$  b(\gamma)$, \  in
\cite{ooz}, Theorem 3.7, the lower bound \vspace{-0.2cm}
      $$ b(\gamma) \geq l_R(\ded/\gamma)(r-1) \hspace{0.5cm} \vspace{-0.2cm} $$
    and in
         \cite{dd}, Proposition 2.1, the upper 
bound   \vspace{-0.2cm} $$  \vspace{-0.2cm} 
b(\gamma) \leq
(r-1) [l_R(R/\gamma)-1] $$ are established.
Hence results $B_1$ and $B_2$ may be viewed as an extension of these bounds
to any ideal $I$.\\  \vspace{-0.2cm}

\noindent There are few cases in which $b(\gamma)\leq r$ \ (see
\ref{coronoto}). A general structure
theorem  for rings satisfying the equality
$b(\gamma)=0$ or $b(\gamma)=1$ \ is presented in \cite{bh}:
    these rings are
called rings {\it
of maximal} \ or {\it   almost maximal length},
respectively.  Note that for $I=\gamma$  the above condition (VC)
becomes: \vspace{-0.2cm}
$$b(\gamma)=0 \  \iff \  r_i=r \   \ for\ each \
\  i=1,....,n. \vspace{-0.2cm}$$  Indeed, the rings of maximal
length   are exactly those
having constant
type sequence. \\
In a series of
recent papers (see
         \cite{dm}, \cite{d}, \cite{dlm}) the authors
attack the problem of classifying rings according to the
value of the
quantity $b(\gamma)$. In the last section we show how type sequences
are an useful
instrument from this point of view, by obtaining a complete classification of
all possible rings having $b(\gamma)\leq r$.

\section{Preliminaries  and notations.}
Throughout this paper    ($R, \m$)   denotes a one-dimensional local Noetherian
      domain with residue field $k$. For 
simplicity, we assume that $k$ is an infinite 
field.
Let $\ \R\ $ be the
integral closure of
$R$ in its quotient field
$K $; we suppose that $\R$ is  a finite $R$-module and a DVR with a
uniformizing parameter $t$,
which means that $R$ is analytically irreducible.  We also
suppose $R$ to be residually rational, i.e., \ $k
\simeq \R/t\R.$ \ We denote the usual    valuation
associated to $\R $
by \vspace{0.1cm} \vspace{-0.2cm}
$$ v : K\longrightarrow \integ \cup
\infty,\qquad v(t)=1 . \vspace{-0.2cm}$$
\noindent In particular \  $v(R)  :=   \{v(a),\;
a\in R,\; a\neq 0\} \subseteq \nat $ \  is the {\it   numerical semigroup
} of $R$.
Under our hypotheses, for any fractional ideals   $I
\supseteq J \neq (0)$   the
length of the $R$-module $I/J$ can be computed by means of valuations  (see
\cite{ms}, Prop. 1):    \vspace{-0.2cm}
  $$l_R(I/J)=\# (v(I) \setminus v(J)).\vspace{-0.2cm}$$

\noindent Given two fractional ideals \ $I, J$ \
we define\   $I:J=\{x\in K \ | \ xJ
\subseteq I\}$.

\begin{(*)}\label{modcan}   \ \\{\rm  In our hypotheses $R$ has
a  {\em canonical module}  \  $\w$,    unique up
to isomorphism. Once for all we  assume   that
$$R\subseteq
\omega \subset \R $$
We shall use the following   properties
(see \cite{hk}):
\begin{enumerate}

\item[(1)] $\w:\w=R$ \ and \  $\omega:(\omega:I)=I$ \ for every fractional
ideal  $I $.

\item[(2)]  $l_R(I/J) = l_R(\w:J/\w:I) $ \ for every fractional
ideals  \   $I \supseteq  J$.

\item[(3)] $R$ is Gorenstein \  {\it if and only if} \ $    \w  =  R$ \
{\it  if and only if} \  $ \ded = R.$   \\
Otherwise\qquad $ \gamma_R \subseteq \ded \subseteq \m.$

\item[(4)] $ v(\w) = \{ j\in \integ \enskip | \enskip c -1
      -j \notin v(R) \}  $, \ hence \ $c-1 \notin v(\w)$ \ and \ $c+\nat
\subseteq v(\w)$.
\item[(5)]   (see \cite{ooz}, Lemma 2.3). For every  fractional ideal $I$, \\
\vspace{-0.3cm}
$$\vspace{-0.2cm} s\in v(I\w)  \quad  if \ and\ only\ if \quad
c-1-s\notin v(R:I).$$
\end{enumerate}}
\end{(*)}

\begin{(*)} \label{ts}  \ \\
{\rm The notion of  {\em type sequence} has been introduced
   by Matsuoka in 1971 and recently revisited in \cite{bdf}; we recall its
definition.
Let $n:=l_R(R/\gamma ) $ \ and let \
\vspace{-0.2cm} $$s_0=0<s_1<\ldots <s_n=c<s_{n+1}=c+1<...\vspace{-0.2cm} $$
be the
      elements of $v(R).$ \ For each $i\geq 1$, \ define the ideal
\vspace{-0.2cm}
  $$R_i:=\{ a\in R\; |\;
v(a)\geq s_i\}. \vspace{-0.2cm}$$
  The chain \ \ \ $R=R_0\supset R_1=\m\supset R_2\supset
\ldots \supset R_n =\gamma_R \supset R_{n+1}\supset... $ \\ induces  the
chain of duals
\vspace{-0.2cm} $$R\; \subset
\; R\; :R_1\; \subset \; .... \subset R\;
:R_n=\; \R \; \subset R\;
:R_{n+1}=t^{-1}\R\subset...  \vspace{-0.2cm}$$
      For every $i\geq 1,$ \ we put \ \ \vspace{-0.2cm}
       $$r_i:=l_R(R:R_{i}/R:R_{i-1})=l_R(\w R_{i-1}/\w R_{i}) \vspace{-0.2cm}$$
and we call \  {\em type\
sequence} \  of $R$
the sequence
$\; [r_1,\dots,r_n] $. \  \vspace{ 0.2cm}\\
      We need in the sequel the
following facts
(see
\cite{bdf}): \begin{enumerate}

\vspace{-0.1cm} \item[(1)] $r:=r_1$ is the {\em Cohen-Macaulay type} of $R$.

\vspace{-0.1cm} \item[(2)]  $ 1\leq r_i \leq  r_1 $ \ for every $\ i\geq 1$ .

\vspace{-0.1cm} \item[(3)] $  \delta   = \sum_1^n \  r_i $,  \ \ and \ \ \ $2
\delta -c = l_R(\w/R)=
\sum_1^n \
(r_i-1)$.
\vspace{-0.1cm}\item[] It follows immediately  that  \  $r-1\leq 2\delta-c$.

\vspace{-0.1cm} \item[(4)] If \ $s_i \in v(\ded),$ then
    $\ r_{i+1}=1 $ (see \cite{ooz}, Prop.3.4).
\vspace{-0.1cm} \item[(5)]$  r_i=1\ $ for every $\ \ i>n$.
\end{enumerate}
}\end{(*)}

\begin{(*)} \ \\
{\rm   A ring $R$ is called  {\it almost Gorenstein} if it
satisfies the equivalent conditions\par
(1) $ \ \m =   \m\ \!\w.$ \par
(2) $\ r -1= 2 \delta - c. $  \\
By the above property \ref{ts},(3), it is clear that $R$ is   {\it almost
Gorenstein} if and only if the type sequence is $
[r,1,\dots,1] $  and that {\em Gorenstein}  means  {\em almost
Gorenstein} with $r=1 $. \\
A ring $R$ is called {\it of maximal length} if it
satisfies the equivalent conditions\par
(1) $ r(c-\delta)=\delta.$ \par
(2) the type sequence is constant $
[r,r,\dots,r] $.   }
\end{(*)}

\begin{(*)}\label{refl}
\ \\{ \rm
      For any proper ideal $I$ of $R$, we denote by \  $\Ib := I\R \cap R$ \
the  \ {\it integral closure}     of    $I$.   Easily we can see that
  \vspace{-0.2cm}
$$ \vspace{-0.2cm}  I \subseteq   I^{**}  \subseteq \w I = \w I^{**}.$$
In fact,    $I^{**}\!\! =
R:\!(R:I) \subseteq \w :(R:I) =  \w I $   \
   and  \  $l_R(\w I^{**}/ \w
I)=l_R(I^{*}/I^{***}\!)=0$. \ Hence \ $ I^{**}  \subseteq \Ib$ \ and \   $
e(  I^{**}) =e(I).$  \ We note also that the
condition
\ $ \w\subseteq I:I  $, \  i.e. \ $ \w I= I,$ \
implies that $ I = I^{**}$. }
\end{(*)}

\begin{(*)}\label{ideals}    \  \\{
\rm
       For any  fractional ideal $I$ we
denote by  $\gamma_I$ the biggest $\R$-ideal
contained in $I$ and    by $\ c_I\ $ the multiplicity of  $\gamma_I$.
        Namely:  \vspace{-0.2cm}
$$\gamma_I:=t^{c_I}\R\subseteq I\
\ {\rm  with} \ \  \
c_I-1\notin v(I),\quad
     R:\gamma_I= t^{c-c_I}\R,\quad v(R:\gamma_I)=\integ_{\geq c-c_I}.  $$
   \noindent  Assume now that $I \subseteq R$ and 
let $\ n_I:=l_R(R/\gamma_I)= c_I-\delta\geq n $. 
\  Then
\begin{enumerate}
\item[(1)]  \ $\gamma_I \subseteq \gamma$ \
and the inclusion \ $  \gamma \subseteq I$  \ implies that \
$\gamma_I =
\gamma$.
\item[(2)]  \ $\sum _{i=1}^{n_I} r_i=  l_R(R: 
\gamma{_I}/R )= c_I-c+\delta$\quad
  and
\item[] \   $\sum_{h=1}^{n_I}   (r_h-1 )=2  \delta-c.$

\item[(3)]  From the square $$\matrix{R&\subseteq&\R\cr | \cap&\ & | \cap
\cr I^*&\subseteq&R:\gamma_I}$$ and the above item we
get
$$l_R(I^*/R)= \sum _{i=1}^{n_I} r_i -l_R(R: \gamma{_I}/I^*)$$

  \end{enumerate}   }
\end{(*)}

\section
{Invariants a(I) and b(I).}
{\rm For any proper ideal $I$ of $R,$ \ we define the two invariants
\begin{enumerate}
\item[] $a(I):=l_R(I^*/R)-l_R(R/I)$
\item[] $b(I):=r l_R(R/I)-l_R(I^*/R)$,
\end{enumerate}
in particular:  $ a(\gamma)= 2 \delta -c ,  \quad 
b(\gamma)=r(c-\delta)-\delta,\quad
    a(\m)= r-1,\quad b(\m)=0.$ \\ \vspace{-0.2cm}

  The aim of the section is to express these invariants   in
terms of the type sequence of $R$. The particular
description given in Theorem \ref{r.nelleformule} allows us to get bounds
and vanishing conditions, improving results of
several authors.  \\ \vspace{-0.2cm}

First we collect    some remarks concerning
   $a(I)$ and $b(I)$.
\begin{rem}\label{remsua-b}  {\rm  Let   $  I $ be a proper ideal of $R$. Then:
\begin{enumerate}
\item[(1)] $a(I)+b(I)=(r-1)\ \! l_R(R/I)$.

\item[(2)]   $a(I)= a(\gamma)-l_R(\w I/I) \leq a(\gamma)$.\par This   easy
computation yields immediately that:\par
   (a) \ \   $a(I)= 0$ \ for every ideal $I$ \ $ \iff R$ is
    Gorenstein\par
(b) \ \   $a(\m)= a(\gamma) \iff R$ is
almost Gorenstein\par
(c) \ \  $I$ {\it canonical}, i.e. $I \simeq \w,
\ \I a(I)= a(\gamma)-l_R(R/\ded)$.\par
\ \ \ \ \ \ \  For a discussion about the invariant \
$\sigma:=a(\gamma)-l_R(R/\ded)$ \ see \par \ \ \ \ \ \ \ \cite{ooz},    3.5,
where we found
examples with \ $\sigma <0$.

\item[(3)]   \
$  b(I)
\geq0 $.  \par This fact   follows by
applying with \ $M=N=R$ \  the J\"ager's
inequality:\vspace{-0.2cm}
$$ \vspace{-0.2cm} l_R(M:I/M:N) \leq l_R(M:\m/M) l_R(N/I)$$
which holds for every fractional ideals $M,\  N,\  I,$ \ such that
$I\subseteq N$ (see \cite{j}, Satz 2).

\item[(4)]   If \ $J \subseteq I, $ \ we have: \par

(a) \ $a(J)-a(I)=l_R(J^*/I^*)-l_R(I/J) $. \par
(b) \ $b(J)-b(I)= r l_R(I/J)- l_R(J^*/I^*)\geq 0.$\par
Assertion (a) is easy to check and (b) follows directly from (a) by means
of (1).
The positivity of  $b(J)-b(I)$ is again a consequence of the Jager's result.
We note in
particular that:
\begin{enumerate} \item[(c)] \  $a(I)= a(I^{**})-l_R(I^{**}/I)$.

\item[(d)] \   $ b(I) =0$ \   for every ideal $I$ containing $\gamma$  if
and only
if     $R$ is a ring
{\it of maximal length}.
\end{enumerate}

\item[(5)]
By  definition \ \
   $\sum_{h=1}^i r_h= l_R(R:R_i/R).$ \ \   Then:

\begin{enumerate}
        \item  $a(R_i)=\sum_{h=1}^i (r_h-1),$ \   in particular   \
$a(R_i)= 2\delta-c $\ \ for every  \
$i\geq n.$
\item   $ b(R_i)=\sum_{h=1}^i
(r-r_h ),$ \ \  in particular \item[]
   for    $i\geq  n$, \ i.e. \
$R_i= t^{c+p}\R,  \ \ p\geq 0$, \ we get \  $b(R_i)= b(\gamma)+ p(r-1).$
\end{enumerate}

\item[(6)]  If $R$ is Arf, i.e. \ $l_R(R:R_i/R)=s_i-i $ \  for every \ $1
\leq i
\leq n$ \ (see \cite{dd}, Proposition 1.15), \   then
$$a(I) \leq (r-1) l_R(R/I) -(i_0 s_1 -s_{i_0}) $$
where $ s_{i_0} $ is the multiplicity of  $I$. \par
In fact, the   hypothesis  $R$  Arf implies that \ $a(R_i)=s_i-2i, \ \ \
b(R_i)= is_1-s_i $. \ Applying  the second formula of (4) to the ideals \
$I \subseteq R_{i_0} =\Ib $
\  we obtain \ $b(I)
\geq i_0 s_1 -s_{i_0}$, \ hence the thesis by (1).

\end{enumerate}
        } \end{rem}

We introduce now another notation.
\begin{notat} \ We associate to any proper ideal 
$I$ the numerical set $V^I$ depending on the 
valuations of
$I^{**}$ \vspace{-0.2cm}
$$\vspace{-0.1cm}V^I:= \{h+1\ | \ h\in \nat \ and   \
s_{h}\in v(I^{**})\}. \vspace{-0.2cm}$$  \end{notat}
The $r_is$ of the {\it type sequence}, with $i\in V^I,$    will be  useful
in our computations.
   \begin{rem}\label{lemma} Let, as usual,\   \  $n_I =
c_I-\delta$.  Then:    \vspace{-0.2cm} $$\# V^I_{\leq n_I}=  l_R(I^{**}
/\gamma_I )  \ \quad and \quad \  \# V^I_{\leq n}=
l_R(I^{**}+\gamma/\gamma ).   \vspace{-0.2cm}  $$
\end{rem}
The basic idea for the next theorem comes from
2.1.(5), which establishes
    a   duality between the
valuations of $\w I$ and those of $I^*\!$.

\begin{thm} \label{diseg}    For any
proper ideal $I$ we have:
\begin{enumerate}
\item $   l_R(I^{**}/\gamma_I )\leq
\displaystyle \sum_{h \leq n_I, \  h \in  V^I}\ r_h  \leq
l_R( R:\gamma_I/  I^*). $
\item $ l_R(I^*/R) \leq \displaystyle\sum_{ \ h\notin
V^I}  r_h = l_R(R/I^{**})+ \displaystyle \sum_{
h\leq n,\ h\notin V^I } (r_h-1).$
    \end{enumerate}    \end{thm}

\underline{Proof}.  \
        The proof is substantially the same as in
\cite{oz2}, Proposition 4.2; some changes are
due to the fact that now we don't
assume that $I$ is a reflexive ideal containing $\gamma$.
        \begin{enumerate}
        \item[(1)]   The first inequality is true by \ref{lemma}, since
$r_h \geq 1 \  $  for
each $h$.\\ For the last one
let $h$ be an integer, $\ \ 1 \leq h \leq n_I.$ \ If $x_{h-1}\in
I^{**}$ is such that $ v(x_{h-1})= s_{h-1}<c_I$,
then
by definition  \vspace {-0.1 cm}$$\vspace {-0.1
cm}\hspace{-0.7cm} r_h= l_R (\w R_{h-1}/\w R_h) =
l_R(x_{h-1}\w+\w R_{h}/\w R_{h})=\#
\{v(x_{h-1}\w +\w R_{h})\setminus v(\w R_h)\}.$$
  Since
$  \
v(x_{h-1}\w) \subseteq v(\w I^{**}) = v(\w I)$, \ by virtue of
\ref{modcan},(5)  the
assignement
\
$y
\to c-1-y$
\ defines an injective map

\vspace{-0.4cm}
$$   \bigcup _{h \in V^I_{\leq n_I} }\{ v(x_{h-1}\w
+\w R_{h})\setminus v(\w R_{h})
\} \longrightarrow \integ_{\geq c-c_I} \setminus v(I^*).$$

\vspace{-0.2cm}
        The conclusion
\ \ $\sum_{h\in V^I_{\leq n_I} }\ r_h \leq l_R(
R:\gamma_I/I^*)$ \ follows, because the sets
\\ $\{ v(x_{h-1}\w +\w R_{h})\setminus v(\w
R_{h})
\}, \ h\in V^I_{\leq n_I} , $ \ are disjoint by construction and because
\ $\integ_{\geq
c-c_I}= v(R:\gamma_I)$.

  \item[(2)] The last inequality in (1) combined with
\ref{ideals} (3) gives:\par
        $l_R(I^*/R)\leq \sum _{h=1}^{n_I} r_h -\sum_{h\in
   V^I_{ \leq n_I}}
r_h=\sum_{h\notin
V^I}r_h=$ \par
\ \ \ \ \ \ \ \ \ \ \ \ \   $=l_R(R/I^{**})+\sum_{h \notin V^I}  (r_h-1 )$.\par
The thesis is now immediate since  $ r_h=1$ \ for
all \ $ h>n $. \ \ $\diamond$

   \end{enumerate}
\begin{coro}\label{aew} \  For any proper ideal
$I$ we have:
\begin{enumerate} \item [] $l_R(\w  I/I)\geq
\displaystyle\sum_{h\in V^I}  (r_h-1)$
\end{enumerate}
\end{coro}
\underline{Proof}. \   By
    \ref{remsua-b}, (4)  and  part (2) of the
theorem,  we obtain
\vspace{-0.2cm} $$\vspace{-0.2cm} a(I)\leq a(I^{**})\leq
\displaystyle \sum_{ \ h\notin V^I  }  (r_h-1) .\vspace{-0.2cm}$$
   Using    \ref{remsua-b}, (3), we conclude   that:
\vspace{-0.2cm}
    $$\vspace{-0.2cm}l_R(\w I/I)= 2\delta-c -
a(I)\geq\sum_{h=1 }^n (r_h-1)- \sum_{h\notin
V^I }   (r_h-1),\vspace{-0.2cm}$$
which is the thesis. \ \ $\diamond$ \\ \vspace{-0.2cm}

\noindent The  last inequality in
Theorem \ref{diseg}, (1)   leads  to
      introduce the following non-negative  invariant.

\begin{defin} \label{def.d} \ For any
proper ideal $I$  we define
\begin{enumerate}
\item[] $d(I):= l_R(R:\gamma_I/I^*) - \displaystyle\sum_{ h\leq n_ I, \ h\in
V^I} r_h$.
\end{enumerate}
        \end{defin}

        It is clear that: \begin{enumerate}\item \
$d(I)=d(uI)$ \ for every unit \ $u \in
\R$; \item \  $d(I)\geq0$, \  by \ref{diseg}; \item \
$l_R(R:\gamma_I/I^*)-rl_R( I^{**}/\gamma_I)\leq d(I)\leq
l_R(R:\gamma_I/I^*)- l_R( I^{**}/\gamma_I)$ \par
  and  the minimal value    is achieved in  a ring   of maximal
length.\end{enumerate}

\begin{coro}\label{dual} Let I be a proper ideal. Then \begin{enumerate}
\item  $ l_R(I/\gamma_I)\leq
l_R( R:\gamma_I/  I^*).$    \item   Equality holds in (1)  $\II\ I$ is
reflexive, \ \ d(I)=0,\ \ \ $r_h=1 \ \ \ \forall \ h\in
   V^I_{\leq n_I} $.\end{enumerate}\end{coro}

\begin{prop} \label{dopodefd}   For any proper ideal $I$
we have:

\begin{enumerate}
\item  \ $d(I)=l_R(\w I/I^{**} )-   \displaystyle{
\sum_{h\in V^I }\ (
r_h-1) }.$
\item
\  $d(I^{**})=d(I)$.

\item \ If \ $ I\subseteq \ded$, \ then \ $d(I)=l_R(\w I/I^{**} )$.

\item \ If \ $ \w\subseteq I:I$,  \ then \ $d(I)=0$.
\item  Let   $\ i_o\in\nat$  be the integer
such that \
$e(I)=s_{i_0}.$ \  Then\par $ d(I)= \displaystyle \sum_{
  h >i_0  ,\ h\notin V^I}
r_h\ -\
l_R(I^*/R_{i_0}^*) $.

\item \ If  $I$ is integrally closed, \ then \ $d(I)=0$.
\item \ If $R$ is almost Gorenstein, then  \ $d(I)=0$.

\end{enumerate}
\end{prop}

\underline{Proof}.
\begin{enumerate}
       \item[(1)]   By (2) of \ref{modcan} \ $l_R( R:\gamma_I/I^*)= l_R(\w I /
\gamma_I)$. \ \ Thus: \par$ d(I)  = l_R(\w I /
\gamma_I)-\displaystyle\sum_{h\in V^I_{\leq n_I}}r_h=
l_R(\w I / I^{**})-(\sum_{h\in V^I_{\leq n_I}}r_h
-l_R(I^{**}/\gamma_I) )=l_R(\w I / I^{**})-
\sum_{h\in V^I }(r_h- 1)  $.
\item[(2)]  It is a consequence of item (1), in view of the fact that \ $\w
I = \w I^{**}$ \ by \ref{refl} and \
$V^I=V^{I^{**}}$
\ by definition.
  \item[(3)] It follows from (1) in view of \ref{ts} (4).
\item[(4)]  The inclusion \ $\w \subseteq I:I $   implies that \ $\w I
=I^{**}$,   hence  the thesis by (3).
\item[ (5)]  After writing \  $l_R(R:\gamma_I/I^*)=
l_R(R:\gamma_I/R_{i_0}^*)-l_R(I^*/R_{i_0}^*),  $     the
thesis is
clear  since \
\vspace{ -0.1cm}
$$l_R(R:\gamma_I/R_{i_0}^*)=\displaystyle  \sum_{
i_0<h\leq n_I}   r_h.
\vspace{ -0.3cm}$$

\item[ (6)] It follows from the above item, because  \ $I=R_{i_0}$.
\item[(7)]  We prove that \ $   \w I = I^{**}$. \ As observed in
\ref{refl}, the inclusion \
$I^{**}  \subseteq  \w I$ \ always holds.
Now \
$\w I  (R:I) \subseteq  \w \m =\m$. \ Thus \ $   \w I \subseteq I^{**}$.
The thesis comes from (1) combined with   the fact
that
$d(I)
\geq0.    \ \ \diamond$

\end{enumerate}

The next theorem extends to any birational overring $S$ of $R$ the
formulas proved in \cite{oz2}
     in the case of the
blowing-up $\LL$ of $R$ along a proper ideal.  We remark also that
    for   $S=\R$ \   the first inequality \ $l_R(S/R)
\leq r \  l_R(R/R:S)$ \
    becomes
        the well-known relation \ $\delta \leq r  (c-
\delta)$.

        \begin{thm} \label{bounds}  \ Let
$S$ be an \ R-overring, \ $R \subseteq S \subseteq
\R$ \ and let
\ $I:=R:S$  \  be its conductor ideal.  Let $i_o\in\nat$ denote the integer
such that \
$e(I)=s_{i_0}$. \
Then:
$$l_R(S/R) =   \displaystyle  \sum_{h\notin
V^I}   r_h-l_R(S^{**}/S)-d(I) \leq
r \ l_R(R/I)\vspace{-0.2cm} $$
$$l_R(S/R) =    \sum_{h\leq i_0}   r_h-l_R(S^{**}/S)+ l_R(S^{**}/R^*_{i_0})
\ \ \ \ \ \ \ \vspace{-0.2cm}$$
        \end{thm}
\underline{Proof}. Since the hypothesis  \ $R \subseteq S \subseteq \R$ \
ensures that \ $\gamma_I=\gamma$, \ the proof of
Theorem 4.4  of
\cite{oz2} works in the general case and we may omit the proof.  \ \
$\diamond$\\ \vspace{-0.2cm}

From Theorem \ref{diseg} we deduce now the following two  formulas which
connect  the invariants
$a(I),\ b(I)$ with the
type-sequence.

\begin{thm}    \label{r.nelleformule}   \
For any proper
ideal $I$ of $R$ we have:
\begin{enumerate}
\item \ $a(I)   = \sum_{ h \notin V^I }\
(r_h-1) - \ l_R(I^{**}/I) \ - d(I)$.
\item \ $b(I)= \sum_{ h \notin V^I}\ (r-r_h) + \  r l_R(I^{**}/I) + d(I)$.
\end{enumerate}
\end{thm}

\underline{Proof}.
\begin{enumerate}
        \item[(1)]  By   \ref{ideals}, (3): \par $a(I)+d(I)+l_R(I^{**}/I)=$ \par
\ \ \ \ \ \ \ \ \ \ \ \ \ \ \
$=l_R(I^*/R)-l_R(R/I)+
        l_R(R:\gamma_I/I^*)-\sum_{h\in V^I_{\leq n_I}} r_h+l_R(I^{**}/I)=$
\par
           \ \ \ \ \ \ \ \ \ \ \ \   $= \sum_{h=1}^{n_I} r_h
-\sum_{h\in
V^I_{\leq n_I}}r_h- l_R(R/I^{**})=$ \par
           \ \ \ \ \ \ \ \ \ \ \ \  $=\sum_{h \notin
V^I}(r_h-1)$.
   \item[(2)]  \ It follows from (1), since $a(I)
+b(I) = (r-1) l_R(R/I)$.
\end{enumerate}

We get immediately interesting lower and upper bounds.

\begin{coro}  \label{coroa}  \ The following inequalities hold:
   \begin{enumerate}
        \item  \  $ a(I) \leq (r-1)l_R( R/ (I^{**}+\ded)) - \
l_R(I^{**}/I).$\par  \ $ a(I) \geq r-1 - \ l_R(I^{**}/I) - d(I).
$
   \item   \ $ b(I) \leq
(r-1)(l_R(R/I)-1)+l_R(I^{**}/I)+d(I).$\par \ $b(I) \geq
(r-1)l_R((I^{**}+\ded)/I) + \ l_R(I^{**}/I).$

\end{enumerate}
\end{coro}
\underline{Proof}.   First recall the positivity of  $d(I) $ and some
properties of type sequences:   \
($i$) $r_h \leq r$ for every $h=1,...,n$; \\ ($ii$) $r_h=1$   for every
$h>n$ and
for every $h$ such that  $s_{h-1} \in v(\ded).$
\\ Then derive   assertions of part    (1) from
the first formula of the theorem. \\
   Since \ $a(I)+b(I)=(r-1)\ \! l_R(R/I)$, \ (2) follows easily from
(1). \ \ $\diamond$ \\ \vspace{-0.2cm}

The first statement in item (1) of the corollary  improves   the inequality
\ $a(I) \leq (r-1) l_R(R/I)   $ \
obtained by
J\"ager in \cite{j}, Korollar 3, (2). \par
The two statements in item (2) generalize  to any ideal $I$ the upper bound  \
  $b(\gamma) \leq (r-1) [l_R(R/\gamma)-1]  $ \ and
  the lower bound
      \ $b(\gamma) \geq l_R(\ded/\gamma)(r-1),  $ \ already known for the
conductor ideal
  (see, respectively,   \cite{dd}, Proposition 2.1 and \cite{ooz}, Theorem
3.7).\par
The second statement in item (1) provides  a  sufficient
condition   for the positivity of $a(I)$. Using \ref{refl} and
\ref{dopodefd},  we have immediately that

\begin{coro}   \label{suff.cond} \ \
   If \ $I$ satisfies the condition \ $\w \subseteq I:I $, \ then \
$a(I)
\geq r-1
\geq 0$.
  \end{coro}

Another direct consequence of \ref{r.nelleformule} is the following.

\begin{coro}   \label{corob}   \
\begin{enumerate}
\item $b(I)\geq r \  l_R(I^{**}/I) \geq 0.$
\item  (Vanishing
condition for  $b(I))$.\vspace{-0.2cm}
$$\vspace{-0.2cm}b(I)=0 \II I=I^{**}, \ r_h=r \   \ \forall   \ \! h \notin
V^I  \ and
\displaystyle
\sum_{h \in V^I\!,\ h\leq n_I } r_h=  l_R(R: \gamma_I/I^*).\diamond $$
\end{enumerate}\end{coro}

Finally we obtain a characterization of the {\it almost Gorenstein}
property in terms of the invariant $a(I)$ (see next $1\II 5)$,  which is
just the analogue  of a theorem   stated by E. Matlis
for Gorenstein rings (see \cite{ma}, Theorem 13.1).

\begin{thm} \label{Matlis13.1}    \ Here "ideal" means "  fractional ideal".
The following facts are equivalent:
\begin{enumerate}
\item  $R$ is almost Gorenstein.
\item   $ \w I =I^{**}$ \ for every non-principal   ideal
$I.$

\item   $l_R(I/J) = l_R(J^*/ I^*)$  \ for every
   reflexive
ideals  $I,J, \ \   J \subseteq I. $
\item    $l_R(I/\gamma_I) = l_R(R: \gamma_I/ I^*)$   \ for
every
reflexive  ideal $I.$
\item    $ a(I) =  (r -1) -  l_R(I^{**}/I)$   \ for
every non-principal
    ideal $I\subseteq R.$
\item   $r-1= 2 \delta -c  $.
\item $\m\w=\m.$
\end{enumerate}
\end{thm}

\underline{Proof}.\begin{enumerate} \item[  (1)]$\hspace{-0.2cm}\I (2)  $ \
As observed in
\ref{refl}, the inclusion \
$I^{**}  \subseteq  \w I$ \ always holds.
Now \
$\w I  (R:I) \subseteq  \w \m =\m$. \ Thus \ $   \w I
\subseteq I^{**}$.

\item[  (2)]$\hspace{-0.2cm}\I (3) $  \  By  2.1:  \ $l_R(J^*/I^*)= l_R(\w
I/\w J)= l_R(I/J)$.

\item[  (3)]$\hspace{-0.2cm} \I (4) $ \ Take \ $J = \gamma_I$.
\item[  (4)]$\hspace{-0.2cm} \I (6) $ \ Take \ $I = \m$.

\item[ (1)]$\hspace{-0.2cm} \I (5) $ \  This implication follows  from
\ref{r.nelleformule},
because in the almost Gorenstein case
$r_h=1$ \   for all \ $h \neq 1$
and \ $d(I)=0$ \ by
\ref{dopodefd}, (7).

\item[ (5)]$\hspace{-0.2cm} \I (6)$ \ \ Take \ $I = \gamma $.
\item[ (6)]$\hspace{-0.2cm} \II (7) \II (1) $  \ These equivalences are
well-known. \end{enumerate}

\section{The special case of $\gamma$.}

The description of the invariant   \ $b:=b(\gamma)$ \  in terms of type
sequence given in Theorem
\ref{r.nelleformule}\vspace{-0.2cm}
$$b =\sum_{h=1}^n(r-r_h)\vspace{-0.2cm}$$
allows us to  complete  the  classification of  all   analitically
irreducible local rings
having
$b\leq r$. Some of the results contained in this section are already
present in the literature (see
  \cite{dm}, \cite{d},  \cite{dlm}). \\ \vspace{-0.2cm}

  From now on we shall denote by
$x\in \m $ an element   such that
$v(x)=e$, \ in other words
$xR$  is a minimal reduction of $\m$.

\begin{lemma}\label{cxR}  \ \\
     Let \ $z:= min\{ y\in
v(R)\ \ |\ \ y\geq c-e\}  $ \ and let \ $  B:= \{h\in
[1,n] \ | \  z<s_h\leq c   \}.$
\begin{enumerate}
\item  $\# B=
l_R(\gamma:_R \m/ \gamma)=l_R (R/\gamma+xR) \geq e-r \geq 1.  $

\item $ \sum_{h\in B} r_h \leq e-1.$
\end{enumerate}
\end{lemma}
\underline{Proof} \ First of all we observe that, called \   $i_0:= min
(B)$,  \ we have by definition \ $z= s_{i_0-1}$ \ and
\ $B=[i_0, n]$.
\begin{enumerate}
\item[(1)]   Obviously we have that\par
\centerline{  $v(\gamma:_R \m)\setminus
v(\gamma)= \{ s_i\in v(R) \ |
\ \ \ c-e\leq s_i<c\}$. }
   Clearly this set is in 1-1 correspondence  with the set
     \par \centerline{$\{i  \ |
    \ z\leq s_i<c\} = [i_0-1,n-1]$, }\ so the first
assertion of (1) is proved.\par
   It is easy to check that  $\qquad x(\gamma:_R \m) = xR\cap \gamma$.\par
     \par
Hence\quad  $l_R \big(xR/x (\gamma:_R \m)\big)=l_R (xR/xR \cap \gamma)=
l_R (\gamma +xR/ \gamma)$ and   to prove the
second equality it suffices to consider the
following  inclusions
$$\matrix{\gamma:_R \m&\subseteq&R\cr \cup | &\ & \cup  | \cr \gamma
&\subseteq&\gamma +xR}$$\par
\par
Finally, since \ $(\gamma +xR)\m \subseteq xR,$ \
we obtain  \ $(\gamma +x\m)  \subseteq xR:\m,$ \
hence \par
\ \ \ \ \ \ $ l_R(\gamma+xR/xR)\leq r $ \  and \par \ \ \ \ \ \  $
l_R(R/\gamma+xR) =
l_R( R/xR)-l_R(\gamma+xR/xR)\geq e-r.$

  \item[(2)] Since\ $c-1 \notin
v(\w )$ \ by (4) of \ref{modcan}, \ \
$v(\w R_{i_0-1})_{<c} \subseteq [c-e, c-2]$. \ Thus:\par
    $ \sum_{h\in B} r_h =l_R (\w R_{i_0-1}/ \gamma) \leq e-1.$ \ \ $\diamond$
\end{enumerate}

\begin{thm}\label{dopoD} \ Let
   $A:=
\{ 1,...,n\}\setminus B.$  \
   The following inequalities hold:\begin{enumerate}\item
$b +e-1 \geq b + \sum_{h\in B}
    r_h  =\sum_{h\in A}
(r-r_h )+r l_R(R/\gamma+xR).$
    \item $b   \geq (r-1)(e-r-1)+\sum_{h\in A}
(r-r_h )$.\end{enumerate}   \end{thm}

\underline{Proof}.
\begin{enumerate}

\item[(1)] We use the description of $b$ in terms of type sequence given in
\ref{r.nelleformule}. \par
$b =\sum_{h=1}^n(r-r_h)=    \sum_{h\in A} (r-r_h )+\sum_{h\in B}(r-r_h )=$
   \par $ \ \  =\sum_{h\in A} (r-r_h
)+r l_R(R/\gamma+xR) -\sum_{h\in B} r_h.$
\item[(2)] Since \
     $ l_R(R/\gamma+xR) \geq e-r,$ \
    by substituting in item (1)  we get\par
$b  \geq \sum_{h\in A}
(r-r_h )+r (e-r)-(e- 1)$, \ which is our thesis. \ \ $\diamond$
\end{enumerate}

\begin{notat}  \ {\rm We denote by}
\begin{enumerate}\item[]  -\ \ $p  $ \ \ the
integer such that \ \ \ $c-e\leq
pe<c$  \ \ $(p \geq 1),$\item []  {\rm and by $g$
the number of gaps of $v(R)$ in the interval
}$(pe, c):$  \item[]
-\ \  $  g = \# \  \nat_{\geq pe} \! \!\setminus\!  v(R),  \quad  (1\leq g
\leq e-1)$.\end{enumerate}
\end{notat}

Formula 1 of Theorem \ref{dopoD} involves the length $  l_R(R/\gamma+xR)$.
  For   the proof of
Theorem \ref{coronoto} we need next two lemmas, which
describe  in detail  the   cases \ $ l_R(R/\gamma+xR)=1, 2 $.

\begin{lemma}\label{delta1}
   The following facts are equivalent:\begin{enumerate}
\item   \  $  l_R(R/\gamma+xR) =1.$
  \item  \ $v(R)=\{0,e,..,pe,c\rightarrow\}.$
\item \ $ts(R)=[e-1,....,e-1,r_n]$.\end{enumerate}\
If $R$ satisfies these equivalent conditions, then  $R$ is a
quasi-homogeneous
singularity  with  \par   $  \delta=c-p-1, \ \ \ b=e(p+1)-c\leq r-1,\ \ \
r=e-1, \ \ \   r_n=e-1-b. $
\end{lemma}

\underline{Proof}. \ $(1)  \II   (2)$ \ is immediate, and  also the fact
that  $R$ is a quasi homogeneous singularity, with
$r=e-1$ by (1) of \ref{cxR}.
\  To prove $(2)\I(3)$,  note that \vspace{-0.3cm}
  $$\vspace{-0.3cm} \sum_{h=1}^{n-1}r_h= l_R(R:R_{n-1}/R)=l_R(x^{-p}R
\cap \R/R)=ep-p= r(n-1)$$
Hence $r_h=r$ \ for each \
$h \in [1, n-1]$. \ Since \ $b= \sum_{h=1}^n(r-r_h)$ \ we get \
$r_n=r-b$. \ Therefore, \
$b<r$ \ and \ $ts(R)=[e-1,....,e-1,e-1-b]$.\\
$(3)\I (2)$ follows, since for each \
$h \in [1, n-1]$ \ the hypothesis \ $r_h=e-1$ \ implies that   \ $
s_h=he$ \
  (see \cite{ooz}, Proposition 4.9).

\begin{lemma} \label{delta}
  \ Assume that  \ $  l_R(R/\gamma+xR) =2.$  \ Then \  $e-2\leq r\leq e-1$ \
and
     there  are   two    possibilities   for   $v(R):$

\begin{enumerate}
\item [(A)]
$   v(R)\!=\! \{0, e,2e,...,ke,y,(k\!+\! 1)e,
y\!+\! e,....,(p\! -\! 1)e, y\! +\! (p\! -\! k\! -\! 1)e, pe,c,\rightarrow\}$
   \item[]   \ \ \ \ \ with  \ \ $  p>k\geq 1,   \
c\leq (p+1)e,$ \   \  $y+(p-k)e\geq
c,$
\ \     $ c-\delta= 2p+1-k $.
\item [(B)] $ v(R)= \{0, e,2e,...,ke,y,(k+1)e,
y+e,....,pe, y+(p-k)e, c,\rightarrow\}$ \item[]
\qquad with \ \ $p\geq k \geq 1, \
c\leq (p+1)e$,    \
$y+(p-k)e<c,$ \ \ \ $ c-\delta= 2p+2-k$.
\end{enumerate}

In both cases we have:

\begin{enumerate} \item[] $\delta= p(e-1)-(p-k )+g=p(e-2)+k+g,$
   \item[] $b+g=r(c-\delta) -p(e-2)-k $ \ and \ $1 \leq g \leq e-2$.
\end{enumerate}

  Moreover:

\begin{enumerate}
   \item  If  \   $r=e-1$, \  then \   $b\geq r+1$ \  and \ $\left
[\begin{array}{ll} case\ (A)& b+g=(p-k+1)e-1 \\ case\ (B)&
b+g=(p-k+2)e-2\end{array}\right.$
\item If  \  \  $r=e-2$, \  then \ $p\geq 2k-1$  \item[] and \
   $\left [ \begin{array}{ll}case\ (A)& b+g=(p-k+1)(e-2)-k \geq k(e-3) \\
case\ (B)& b+g=(p-k+2)(e-2)-k > k(e-3).
\end{array}\right.$
\end{enumerate}

\end{lemma}

\underline{Proof}. \  The fact that \  $e-2\leq r \leq e-1 $ \ follows
immediately from    \ref{cxR}.(1).
\begin{enumerate}
\item[(1)] In case
(A) \par $b+g=(e-1)( 2p+1-k) -p(e-2)-k= (p-k+1)e-1$. \par Then the
inequality \
$g  \leq e-2$ \ leads to \ $b \geq r+2$.\par   In case (B) \par  $b+g=(e-1)(
2p+2-k) -p(e-2)-k=(p-k+2)e-2$ \par
   and the same inequality leads to \ $b \geq r+1$.

\item[(2)]  It suffices to prove that \
$2y<c+e$;
\ in fact from this we can deduce  that \vspace{-0.2cm}
$$ 2ke<2y<c+e \leq (p+2) e, {\rm \ hence  \ } \ 
p>2k-2.\vspace{-0.2cm}$$\ If \ $2y
\geq c+e$, \ then by considering the structure of $v(R)$ we can easily see
that \ $\m^2 \subseteq t^e \m$. \ Thus, \ $\m=t^e(R:\m)
\subseteq R
\subseteq R:
\m$,
\ contradicting the assumption \ $r=e-2$. \ \ $\diamond$
\end{enumerate}

\begin{coro} \label{aggiunto} \   Assume that  \ 
$b<q(r-1),\ \ q\geq 1$, \ then  \vspace{-0.2cm}
$$e-r\leq \ l_R(R/\gamma+xR)\leq q.   \vspace{-0.2cm}$$
  In particular \begin{enumerate}
\item  $\ \  \ 0\ \ \leq b <\ r-1\  \I \ r=e-1\ $ \ and \ $l_R(R/\gamma+xR)=1.$
\item
$ \ r-1< b<2r-2 \I  e-2\leq r\leq e-1  $ \ and \ $
l_R(R/\gamma+xR)=2.$\end{enumerate}
\end{coro}
\underline{Proof}.   \ Item (2) of \ref{dopoD} implies that
$(r-1)(e-r-1-q)<0$,\ so $e-1-q<r   $  and item (1) gives \ $
rl_R(R/\gamma+xR)< e-1+q(r-1)<r(q+1);$ \ hence the thesis using also
\ref{cxR} (1).  \\
(a) is the case $q=1$, (b) is the case $q=2$, with the further assumtion
$b>r-1$. It suffices to recall that  by \ref{delta1} \
$ l_R(R/\gamma+xR)=1\I b\leq r-1$.\\ \vspace{-0.2cm}

From these technical observations and Theorem \ref{dopoD} we   deduce  the
statements   of the next theorem,  which are
partially already known (see  \cite{bh}, \cite{d}, \cite{dlm},
\cite{dm}). Nevertheless, they give a    complete   classification of  all
analitically irreducible local rings
having
$b\leq r$. \par
We shall
consider separately the   cases: \ $1) \   b<
r-1;
   \ 2) \ b=r-1; \ 3) \ b=r.$

\begin{thm}
\label{coronoto} Suppose $R$ not Gorenstein.
\begin{enumerate}
\item The following facts are equivalent:\par
   (a) $ b< r-1 $\par
(b) $v(R)=\{0,e,..,pe,c\rightarrow\}$ \ with \ $pe+2<c \leq (p+1)e  $ \par
  (c) $ts(R)=[ e-1, e-1,...,e-1, r_n], \ r_n>1  $. \par
If these conditions  hold, then
\item[]   $  l_R(R/\gamma+xR) =1,\   \
c=(p+1)e-b,  \ \ r=e\!-\!1,\  r_n=e \!-\!1\!-\! b.$

\item $b =r-1\I  \left \{ \begin{array}{l}
r=e-1    \\ or\\
   r=e-2 \end{array}\right.$

$1^{st}$ case)  \ \ The following facts are equivalent:\par
   \ \ \ \ (a) $ b= r-1=e-2 $\par
\ \ \ \ (b) $v(R)=\{0,e,..,pe,c\rightarrow\}$ \ with \ $ c=pe+2  $ \par
  \ \ \ \ (c) $ts(R)=[ e-1, e-1,...,e-1, 1]$. \par
\item[] \ \ \ \ If these conditions  hold, then \ $ l_R(R/\gamma+xR)= 1. $

$2^{nd}$ case) \ \ The following facts are equivalent:\par
   \ \ \ \ (d) $ b= r-1=e-3 $\par
\ \ \ \ (e)  either \ $v(R)=\{0,e,2e-1,2e,3e-1\rightarrow\} $  \par
\ \ \ \ \ \ \ \ \ or \ \ \ \ \ $v(R)=\{0,e, y,2e\rightarrow\}$
\ with \ $ e<y\leq e+ \frac {e-1}{2}$ \par
  \ \ \ \ (f) either \ $ts(R)=[e-2,e-2,r_3,r_4]$ \ with \ $r_3+r_4=e-1$ \par
\ \ \ \ \ \ \ \ \ or \ \ \ \ \ $ts(R)=[ e-2,r_{2},r_3]$ \  with
\
$r_{2}\!+r_3\!=e-1$.
\item[] \ \ \ \ If these conditions  hold, then \ $ l_R(R/\gamma+xR)= 2. $

\item $b =r \I  \left \{ \begin{array}{l}
(g) \ \ \  r=e-2, \ l_R(R/\gamma+xR) =2 \\ or\\
(j) \ \ \ r=2,\   e=5, \ l_R(R/\gamma+xR) =3 .\end{array}\right.$

In   case  $(g),\ \ v(R)$ is one of the  following  sets\par
$  \{0,4,8,9,12,13,16\rightarrow\};
\\    \{0,4,8,11,12,15,16,19\rightarrow\};
\\   \{0,e, 2e-2, 2e, 3e-2\rightarrow\} ,\ \ with\ \ e\geq 4;
\\ \{0,e, e+z, 2e-1 \rightarrow\},\ \ with\
\ 0<z\leq  \displaystyle{\frac {e-2}{2}}
, \ \ e\geq 4$.\par

In          case  $(j),\  v(R)$ is one of the  following  sets  (see \cite{dm},
Rem. 2.7)\par
$  \{0,5,6,7,10 \rightarrow\};
\\    \{0,5,6,8,10 \rightarrow\};
\\   \{0,5,8,9,10,13 \rightarrow\}$

\end{enumerate}\end{thm}

\underline{Proof}   \begin{enumerate}
\item[(1)]  $(a) \I (b)$.  If  $b<r-1$, then by \ref{aggiunto},2,
$l_R(R/\gamma+xR)= 1$. \\ By Lemma \ref{delta1} \ $v(R)= \{ 0,
e, 2e,...,   pe,   c, \rightarrow\
   \} $   with $ \  (p+1)e\geq c$. \ Then \ $b=(p+1)e-c<e-2$ \ implies that
\ $pe+2<c$.\par

$(b) \I (c)$. The hypothesis   $pe+2<c$  gives \ $b<r-1$.   Then
$r_n=r-b>1$.\par
$(c) \I (a)$. We have \ $b = r-r_n$, \ hence the thesis.

\item[(2)] By substituting   $b =r-1$ in Formula 2 of    \ref{dopoD}, we
get $(r-1)(e-r-2)\leq 0.$
Two cases are  possible: \ $r=e-1$ \  or \ $r=e-2$ \ and
\ $\sum_{h\in A}(r-r_h)=0$. \par
First case.\par
$(a) \I (b)$.  As in (1) one gets \ $l_R(R/\gamma+xR)= 1$. \ Then \
$v(R)=\{0,e,..,pe,c\rightarrow\}$ \ and \
$b=(p+1)e-c=e-2$, \ hence \ $c=pe+2$.\par
$(b) \I (c)$. See Lemma \ref{delta1}.\par
$(c) \I (a)$. In fact, \ $b=r-r_n=r-1$. \par
Second case. \par
$(d) \I (e)$. If \ $b=r-1$ \ and \ $r=e-2$, \ then by (2) of \ref{dopoD} we
have \ $ \sum_{h\in A} (r-r_h)=0$ \ and from item
(1) of
\ref{dopoD} we obtain \ $ l_R(R/\gamma+xR) =2$.  \\
  It follows that \
$t.s(R)=[e-2,....,e-2,r_{n-1},r_n]$.\par
We have to consider the two cases of  Lemma  \ref{delta}.\par
In case (A) with $b= e-3,$  \ from the inequality \ $g \geq k(e-3)-b$
we get \vspace{-0.2cm}
   $$\ e-2\geq g \geq
(k-1)(e-3) \vspace{-0.2cm}$$ Three possibilities occur:\par
1) \ $ k=1$. \ Then \ $p=2,\
g=e-2, \ c=3e-1, \  y=2e-1 $. \ In conclusion \ \vspace{-0.2cm}
$$v(R)=\{0,e,2e-1,2e,3e-1\rightarrow\}.\vspace{-0.2cm}$$
2) \  $  k=2$. \ Then \ $p=3, \ g=e-3, \
c= 4e-2, \ y=3e-2$, \ so \ $2y> c+e$, \ absurd (see (2) in the proof of
\ref{delta}). \par
3) \ $  e=4 , \ k=3$. \ Then \ $p=5,
\ g= 2, \  c=23, \ y=15,$ \ as above
   impossible   since
$2y> c+e.$

  In case (B) with  $b= e-3,$ \  since \ $g \geq k(e-3)+e-2-b$, \
we obtain \vspace{-0.2cm}
   $$\ e-2\geq g \geq
  k (e-3)+1 \vspace{-0.2cm}$$
The only possibility is \ $  k=1$. \ Then we get \ $p=1, \ g=e-2  $ \ and \
  $v(R)=\{0,e, y,2e\rightarrow\}$ \  with \
$ e<y\leq e+ (e-1)/2 $.\par
$(e) \I (f)$. Let $R_0$ be the monomial ring such that \
$v(R_0)=v(R)=  \{0,e,2e-1,2e,3e-1\rightarrow\} $. \ Then \ $r(R) \leq
r(R_0) =e-2$. \  Since  \ $l_R(R/\gamma+xR)= 2$, \ we
have by item (2) of \ref{delta} \
  $r(R) \geq e-2$. \  Then  \ $r(R)  =e-2$. \  We easily compute \
$b=(c-\delta)r-\delta=e-3$. \ Substituting in   item (2) of
\ref{dopoD} we obtain  \ $\sum_{h\in A} (r-r_h)=0$, \ hence \ $r_2=e-2$ \
and \
$r_3+r_4=e-1$.  \par
The same reasoning holds for \  $v(R)=\{0,e, y,2e\rightarrow\}$. \par
$(f) \I (d)$. It suffices to recall that \ $b= \sum_{h=1}^n(r-r_h)$.

\item[(3)] Assume \ $b=r$. \  From  item (2) of \ref{dopoD}  it follows
that \ $(r-1)(e-r-2)\leq 1$, \\ then  using also \ref{dopoD},(1) we argue
that either \
$r=2$ \ and \
$e\leq 5$,
\ or
\
$ l_R(R/\gamma+xR) =2 $ \ and \ $  r\geq e-2$.\\
Since the cases \  $r=2, \ e=3$ \ and  \ $l_R(R/\gamma+xR) =2, \  r= e-1$ \
are
impossible by Lemma \ref{delta}, \ the first 
assertion is proved. \\ \vspace{-0.2cm}

Case (g): \ \
  \  $  l_R
(R/\gamma+xR) =2$ \  and \ $ b= r=e-2$.  \\ We proceed   analogously
to the proof of (2).\\
In case (A) we have \vspace{-0.2cm}
$$e-2\geq g= (p-k)(e-2)-k\geq  (k-1)(e-3)-1 \vspace{-0.2cm} $$
This gives the following possibilities:\par
1) \ $k=1$. \ Then \ $p=2\  g=e-3,   \ c=3e-2, \ y=2e-2$.
\ Hence \vspace{-0.2cm} $$ v(R)=\{0, e,2e-2,2e,
3e-2\rightarrow\},\quad e\geq 4. \vspace{-0.2cm}$$
2)  \ $k=2$.\par
\ \ \ i) $  k=2 , \ p=4=e$. \ Then \ $g=2, \
c-\delta= 7, \ \delta= 12,\  c= 19, \ y=11$,  \vspace{-0.2cm} $$\
v(R)= \{0,4,8,11,12,15, 16, 19 \rightarrow\} \vspace{-0.2cm}$$
\ \ \ ii) $  k=2 , \ p=3$. \ Then \ $g=e-4, \ c=4e-3,  \ y \geq 3e-3 \I
2y>c+e$ \par \ \ \ \ \  impossible. \par
3) \  $ k=3$.  \par
\ \ \ i)  \ $ k=3, \ e=5$. \ Then \ $p=5,\ g=3, c= 29, \ y=19  \I 2y>c+e$
\par \ \ \ \ \  impossible. \par
\ \ \ ii) \ $ k=3, \ e=4$. \ Then \  $p=5,\ g=1, c= 22, \ y=14  \I 2y>c+e$,
\par \ \ \ \ \   impossible. \par
4) \ $  k=4,  \ e=4$. \ Then \ $p=7, \  g=2, \ c=31, \ y=19  \I 2y>c+e$
\par \ \ \ \ \  impossible. \par
  In case (B) we have  \vspace{-0.2cm} $$\ e-2\geq
g=(p-k+1)(e-2)-k\geq k(e-3) \vspace{-0.2cm} $$ 
and the following possibilities:\par
1) \ $  k=1$. \ Then $ p=1,  \ g=e-3, \ c=2e-1$. \ Hence \vspace{-0.2cm}
$$ v(R)=\{0,e,
e+z, 2e-1 \rightarrow\}  \ \ {\rm with} \ \ 0<z\leq
{\frac {e-2}{2}}, \ \ e\geq 4 \vspace{-0.2cm}$$
  2) \  $ k=2, \ e=4$. \ Then $ \  p=3,\ \ g=2, \
v(R)=\{0,4,8,11, 12,15,16\rightarrow \}.$  \par
Case (j) is treated in \cite{dm}. \ \ \ $\diamond$
\end{enumerate}


\begin{thebibliography}{99}

\bibitem{bdf} V. Barucci, D. E. Dobbs, M. Fontana,
Maximality Properties in
Numerical
Semigroups
       and Applications to One-Dimensional Analytically
Irreducible Local
Domains,   Mem. Amer.
Math. Soc.   vol. 125, n. 598 (1997).

\bibitem{bf} V. Barucci, R. Fr\"oberg, One-Dimensional
Almost Gorenstein
Rings,   Journal of
Algebra 188  (1997) 418-442.

\bibitem{b} R. Berger, Differentialmoduln eindimensionaler
lokaler Ringe,
Math. Zeitschr. 81,
326-354 (1963).

\bibitem{bh} W. C. Brown, J. Herzog,   One Dimensional Local
       rings of Maximal and Almost Maximal Length,   Journal
of Algebra  151,
332-347 (1992).

\bibitem{dd}  M. D'Anna, D. Delfino, Integrally closed
ideals and type
sequences in
one-dimensional
       local rings,   Rocky Mountain J. Math.  27,  (4)
(1997)  1065-1073.


\bibitem{dm}   M. D'Anna, V. Micale, Construction of one-dimensional rings
with fixed value of \ $
  t(R)
\lambda_R(R/C)-\lambda_R(\R/R)$, International 
Journal of Commutative rings  2 (1) (2002).

\bibitem {d} D. Delfino, On the inequality $\lambda(\R/R)
\leq t(R)
\lambda(R/C)$ for
one-dimensional local rings, J. Algebra 169  (1994),
332-342.

\bibitem {dlm} D. Delfino, L. Leer, R. Muntean, A length
inequality for
one-dimensional local rings, Comm. Algebra 28 (2000),
2555-2564.

\bibitem{hk} J. Herzog, E. Kunz,  Der kanonische Modul
eines Cohen-Macaulay
Rings, Lecture
Notes in Math. vol.  238, Springer, Berlin,  (1971).

\bibitem{j} J. J\"ager, L\"angenberechnung  und kanonische
Ideale in
       eindimensionalen Ringen,   Arch. Math. 29  (1977),
504-512.

\bibitem{ma} E. Matlis, 1-Dimensuional
Cohen-Macaulay Rings, Springer-Verlag (1973).

\bibitem{ms} T. Matsuoka, On the degree of singularity of
one-dimensional
analytically
irreducible noetherian rings, J. Math. Kyoto Univ. 11:3
(1971), 485-491.

\bibitem{oz1} A. Oneto, E. Zatini, \  Type-sequences of
modules, \  J.
Pure  Appl. Algebra, 160 (2001), 105-122.

\bibitem{ooz} F.Odetti, Anna Oneto, Elsa Zatini, \ Dedekind
Different and
Type Sequence,  \ Le Matematiche, Vol.LV (2000),
Fasc.II, 467-484.

\bibitem{oz2} A. Oneto, E. Zatini, \  An application of
type sequences to
the blowing-up, to appear on Beitr\"age zur Algebra
und Geometrie.

\end{thebibliography}
\end{document}